\documentclass[12pt]{amsart}

\usepackage{amssymb}
\usepackage{amsmath}
\evensidemargin\oddsidemargin

\begin{document}

\setcounter{page}{1}

\newtheorem{theorem}{Theorem}

\newtheorem{THEO}{Theorem\!\!}
\newtheorem{Rem}{Remark\!\!}
\newtheorem{Frem}{Final Remark\!\!}
\newtheorem{Corun}{Corollary 1\!\!}
\newtheorem{Proni}{Proposition 2\!\!}
\newtheorem{Prosan}{Proposition 3\!\!}
\newtheorem{ProV}{Proposition 4\!\!}
\newtheorem{Cogo}{Corollary 5\!\!}
\newtheorem{Corroku}{Corollary 6\!\!}
\newtheorem{Coshichi}{Corollary 7\!\!}
\newtheorem{Proachi}{Proposition 8\!\!}
\newtheorem{Prokyu}{Proposition 9\!\!}

\renewcommand{\theTHEO}{}
\renewcommand{\theRem}{}
\renewcommand{\theFrem}{}
\renewcommand{\theCorun}{}
\renewcommand{\theProni}{}
\renewcommand{\theProsan}{}
\renewcommand{\theProV}{}
\renewcommand{\theCogo}{}
\renewcommand{\theCorroku}{}
\renewcommand{\theCoshichi}{}
\renewcommand{\theProachi}{}
\renewcommand{\theProkyu}{}

\newcommand{\eqnsection}{
\renewcommand{\theequation}{\thesection.\arabic{equation}}
    \makeatletter
    \csname  @addtoreset\endcsname{equation}{section}
    \makeatother}
\eqnsection

\def\a{\alpha}
\def\CC{{\mathbb{C}}} 
\def\Ea{E_\a}
\def\Eac{{\mathcal E}}
\def\EE{{\mathbb{E}}} 
\def\elaw{\stackrel{d}{=}}
\def\eps{\varepsilon}
\def\Fa{F_\a}
\def\Ga{G_\a}
\def\hta{{\hat \tau}}
\def\htta{{\hat {\tilde \tau}}}
\def\hT{{\hat T}}
\def\hX{{\hat X}}
\def\ii{{\rm i}}
\def\lbd{\lambda}
\def\lacc{\left\{}
\def\lcr{\left[}
\def\lpa{\left(}
\def\lva{\left|}
\def\NN{{\mathbb{N}}} 
\def\pb{{\mathbb{P}}}
\def\rl{{\mathbb{R}}}
\def\racc{\right\}}
\def\rcr{\right]}
\def\rpa{\right)}
\def\Ta{T_\a}
\def\Ua{U_\a}
\def\Un{{\bf 1}}

\def\d{\, \mathrm{d}}
\def\qed{\hfill$\square$}
\def\elaw{\stackrel{d}{=}}

\newcommand{\fin}{\vspace{-0.3cm}
                  \begin{flushright}
                  \mbox{$\Box$}
                  \end{flushright}
                  \noindent}
                  
\title[Hitting densities for stable processes]
      {Hitting densities for spectrally positive stable processes}

\author[Thomas Simon]{Thomas Simon}

\address{Laboratoire Paul Painlev\'e, U. F. R. de Math\'ematiques, Universit\'e de Lille 1, F-59655 Villeneuve d'Ascq Cedex. {\em Email} : {\tt simon@math.univ-lille1.fr}}



\begin{abstract} A multiplicative identity in law connecting the hitting times
  of completely asymmetric $\a-$stable L\'evy processes in duality is established. In the spectrally positive case, this identity allows with an elementary argument to compute fractional moments and to get series representations for the density. We also prove that the hitting times are unimodal as soon as $\a\le 3/2.$ Analogous results are obtained, in a much simplified manner, for the first passage time across a positive level. 
\end{abstract}

\maketitle

\section{A multiplicative identity in law}

Let $\{X_t, \, t\ge 0\}$ be a spectrally positive L\'evy $\a-$stable process ($1< \a <2$), starting from zero and normalized such that
\begin{equation}
\label{normal}
\EE\lcr e^{-\lbd X_t}\rcr \; =\; e^{t\lbd^\a}, \quad t, \lbd \ge 0.
\end{equation}
For every $x > 0,$ consider the first exit times 
$$T_x \, =\, \inf\{ t > 0, \; X_t > x\}, \quad \hT_x \, =\, \inf\{ t > 0, \; X_t < -x\}$$
and the first hitting times 
$$\tau_x\, =\,  \inf\{ t > 0, \; X_t = x\}, \quad \hta_x \, =\, \inf\{ t > 0, \; X_t = -x\}.$$
Notice that since $X$ has no negative jumps, one has $\hta_x = \hT_x$ a.s. and $\hta_x$ is a positive stable random variable with index $1/\a$ - see e.g. Theorem 46.3 in \cite{S}. Besides, (\ref{normal}) and the optional sampling theorem give the normalization
$$\EE\lcr e^{-\lbd \hta_x}\rcr \; =\; e^{-x\lbd^{1/\a}}, \quad x, \lbd \ge 0.$$
On the other hand, the process $X$ crosses the level $x > 0$ by a jump so that $X_{T_x} > x$ and $\tau_x > T_x$ a.s. Introducing the overshoot $K_x \, =\, X_{T_x} \, -\, x,$
the strong Markov property at time $T_x$ entails that $\{ X^x_t = X_{T_x +t} - X_{T_x}, \; t\ge 0\}$ is independent of $(T_x, K_x).$ In particular, setting $\hta^x_y = \inf\{ t > 0, \; X^x_t = -y\}$ for $y >0$, we have
$$\tau_x \; = \; T_x \; +\; \hta^x_{K_x}\; \elaw\; T_x \; +\; K_x^\a \htta_1\; \elaw\; x^a( T_1 \; +\; K_1^\a \htta_1)$$
where $\htta_1$ is an independent copy of $\hta_1$ and the two identities in law follow at once from the self-similarity relationships 
$$(T_x, \tau_x, \hta_x) \elaw x^\a (T_1, \tau_1, \hta_1)\quad\mbox{and}\quad K_x \elaw x K_1,$$ which are themselves simple consequences of the self-similarity of $X$ with index $1/\a.$ Specifying the above to $x=1$ yields the basic identity 
\begin{equation}
\label{Over}
\tau_1 \; \elaw \; T_1 \; +\; K_1^\a \htta_1,
\end{equation}
which has been known for a long time \cite{Do}. Notice that the laws of the three random variables appearing in the right-hand side of (\ref{Over}) are more or less explicit: the density of $\htta_1$ is that of a positive $1/\a-$stable variable, the density of $K_1^\a$ comes from a classical computation by Port involving some Beta integral - see e.g. Exercise VIII.3 in \cite{B} - and the density of $T_1$ can be obtained in changing the variable in the main result of \cite{BDP} thanks to the identity $T_1 = S_1^{-\a},$ with the notation $S_t = \sup \{X_s, \; s\le t\}$ for the supremum process. However, lack of explicit information on the law of the bivariate random variable $(T_1, K_1)$ prevents from applying (\ref{Over}) directly to derive an expression of the law of $\tau_1$. Recently Peskir \cite{Pe} circumvented this difficulty with an identity of the Chapman-Kolmogorov type linking the laws of $X_1, S_1$ and $\tau_1$, and providing an explicit series representation at $+\infty$ for the density of $\tau_1.$ In this paper we will follow the method of our previous article \cite{Th1} in order to derive, firstly, a simple identity in law for $\tau_1.$ Our main result reads:

\begin{THEO} One has 
\begin{equation}
\label{Main}
\tau_1\; \elaw\; \Ua \times \hta_1,
\end{equation}
where $U$ is an independent, positive random variable with density
$$f_{\Ua}(t)\; =\;\frac{-(\sin \pi\a) t^{1/\a}}{\pi(t^2 - 2t \cos \pi\a  +1)}\cdot$$
\end{THEO}

\noindent
{\em Proof}: From (\ref{Over}) we have for any $\lbd \ge 0$
\begin{eqnarray*}
\EE\lcr e^{-\lbd^\a \tau_1}\rcr \; = \; \EE\lcr e^{-\lbd^\a (T_1 \; +\; K_1^\a \htta_1)}\rcr & = & \EE\lcr e^{-\lbd^\a T_1} \EE\lcr e^{-\lbd^\a K_1^\a \htta_1}\, \vert \, (T_1, K_1)\rcr\rcr\\
& = & \EE\lcr e^{-(\lbd^\a T_1 +\lbd K_1)} \rcr\; =\; \EE\lcr e^{-(T_\lbd + K_\lbd)} \rcr,
\end{eqnarray*}
where in the third equality we recalled that 
$\htta_1$  is a $(1/\a)-$stable positive variable with given normalization. On the other hand, the double Laplace transform of $(T_1, K_1)$ had been evaluated long ago by Fristedt, and yields 
the simple expression
$$\int_0^\infty e^{-s\lbd} \EE\lcr e^{-(T_\lbd + K_\lbd)} \rcr d\lbd\; =\; \lpa\frac{1}{s-1}\rpa - \lpa\frac{\a}{s^\a -1}\rpa$$
for every $s >1$ - see (2.16) in \cite{Do}. It is possible to invert the right-hand side in noticing that 
$$ \frac{1}{s-1}\; =\; \int_0^\infty e^{-s\lbd} e^{-\lbd}
d\lbd\quad \mbox{and}\quad \frac{1}{s^\a-1}\; =\; \int_0^\infty e^{-s\lbd} \lpa \a\lbd^{\a-1} \Ea'(\lbd^\a)\rpa d\lbd$$
for every $s >1,$ where $\Ea'$ is the derivative of the Mittag-Leffler function
$$\Ea(z)\; =\; \sum_{n=0}^{\infty} \frac{z^n}{\Gamma (1 +\a n)}, \quad z\in\CC.$$
Above, the computation involving $\Ea$ had been made originally by Humbert and we refer to the discussion after (5) in \cite{Th1} for more details. Putting everything together gives an expression of the Laplace transform of $\tau_1$:
\begin{equation}
\label{Second}
\EE\lcr e^{-\lbd^\a \tau_1}\rcr \; = \;F_1(\lbd)\; -\; \a\Fa'(\lbd), \quad \lbd \ge 0,
\end{equation}
where we used the notation $\Fa(x) = \Ea(x^\a)$ (notice that $F_1(x) = e^x$) for every $x\ge 0.$ We will now prove that the function $\Ga : \lbd\mapsto F_1(\lbd) - \a\Fa'(\lbd)$ is completely monotone, reasoning as in Theorem 1 in \cite{Th1}. The well-known curvilinear representation of the analytic function $1/\Gamma$ yields
$$F_1(\lbd)\; =\; \frac{1}{2\pi \ii} \int_{H_\lbd} \frac{e^t}{t-\lbd} dt\quad\mbox{and}\quad \a\Fa'(\lbd)\; =\; \frac{1}{2\pi \ii} \int_{H_\lbd} \frac{\a\lbd^{\a-1}e^t}{t^\a-\lbd^\a} dt$$
for every $\lbd >0,$ where $H_\lbd$ is a Hankel path encircling the disk centered at the origin with radius $\lbd$ - see e.g. Chapter IX.4 in \cite{D}. We deduce
$$\Ga(\lbd)\; =\; \frac{1}{2\pi\ii} \int_{H_\lbd}e^t g_\lbd(t) dt, \quad \lbd > 0$$
where 
$$g_{\lbd}(t)\; =\; \lpa \frac{1}{t-\lbd}\rpa\; -\; \lpa \frac{\a\lbd^{\a -1}}{t^\a -\lbd^\a}\rpa$$
is analytic on $\CC / \{\lbd\}\cup(-\infty, 0].$ Besides, since $g_\lbd(t)\to
(\a -1)/2\lbd$ as $t\to\lbd,$ it can be continued on the whole complex plane cut on the negative real axis, which entails
$$\Ga(\lbd)\; =\; \frac{1}{2\pi\ii} \int_{H}e^t g_\lbd(t) dt, \quad \lbd > 0$$
where $H$ is a Hankel path independent of $\lbd.$ We can then compute, as in Theorem 1 in \cite{Th1}, 
\begin{eqnarray*}
\Ga(\lbd) & = & \lim_{\rho\to 0} \frac{1}{2\pi \ii}\lpa  \int_\rho^\infty e^{-s} \lpa \frac{1}{s+\lbd} + \frac{\a\lbd^{\a -1}}{e^{\ii \pi\a}s^\a  - \lbd^\a}\rpa - e^{-s} \lpa \frac{1}{s+\lbd} + \frac{\a\lbd^{\a -1}}{e^{-\ii \pi\a}s^\a  - \lbd^\a}\rpa ds \rpa\\
& = & \frac{-\a\sin\pi\a}{\pi}\int_0^\infty e^{-s} \lpa \frac{\lbd^{\a-1} s^\a}{s^{2\a} - 2s^\a\lbd^\a\cos\pi\a +\lbd^{2\a}}\rpa ds\\
& = & \frac{-\a\sin\pi\a}{\pi}\int_0^\infty e^{-\lbd s} \lpa \frac{s^\a}{s^{2\a} - 2s^\a\cos\pi\a +1}\rpa ds, 
\end{eqnarray*}
so that from (\ref{Second}) and the expression of the density $f_{\Ua}$,
$$\EE\lcr e^{-\lbd^\a \tau_1}\rcr \; = \; \EE\lcr e^{-\lbd \Ua^{1/\a}}\rcr$$
for all $\lbd \ge  0.$ Reasoning now exactly as in Theorem 3 in \cite{Th1} entails finally the desired identity in law $\tau_1\elaw \Ua \times \hta_1.$

\fin

\begin{Rem} {\em In Theorem 1 of \cite{Do}, the double Laplace transform of the function $(t,x) \mapsto \pb[\tau_x \ge t]$ was computed for general L\'evy processes with no negative jumps. Our result can be viewed as an inversion of this Laplace transform in the particular stable case.}

\end{Rem}

\section{Fractional moments and series representations}

In this section we will obtain two expressions for the density of $\tau_1$ (a density which is readily known to exist from our main result, see also \cite{M} for the general context) in terms of convergent series, the first one being useful in the neighbourhood of $0$ and the second in the neighbourhood of $+\infty.$ Finding different asymptotic expansions for densities of first passage times is a classical issue in probability, see Chapter 1 in \cite{Fre}. Recall also that this work has been done long ago for $\hta_1,$ respectively by Pollard and Linnik - see (14.31) and (14.35) in \cite{S}. Our approach is different from the above references and relies on Mellin inversion. We will hence first compute the fractional moments of $\tau_1,$ as a simple consequence of our main result:

\begin{Corun} 
\label{Mom1}
For any $s\in (-1-1/\a, 1-1/\a),$ one has
$$\EE\lcr \tau_1^s\rcr\; =\; \frac{\Gamma(1-\a s) \sin (\pi (\a -1)(s+1/\a))}{\Gamma(1-s) \sin(\pi(s+1/\a))}\cdot$$
\end{Corun}

\noindent
{\em Proof}: From (\ref{Main}) we have $\EE\lcr \tau_1^s\rcr =\EE\lcr \hta_1^s\rcr\EE\lcr \Ua^s\rcr$ for every $s\in\rl.$ The computation of $\EE\lcr \hta_1^s\rcr$ is easy and classical:
\begin{eqnarray*}
\EE\lcr \hta_1^s\rcr & = & \frac{1}{\Gamma(-s)}\int_0^\infty \lbd^{-(s+1)} \EE\lcr e^{-\lbd \hta_1}\rcr d\lbd\; = \; \frac{1}{\Gamma(-s)}\int_0^\infty \lbd^{-(s+1)} e^{-\lbd^{1/\a}} d\lbd\; =\; \frac{\Gamma(1-\a s)}{\Gamma(1-s)}
\end{eqnarray*}
for every $s < 1/\a.$ Similarly, one can show that 
$$\EE\lcr \Ua^s\rcr \; =\; \frac{1}{\Gamma(-\a s)} \int_0^\infty \lbd^{-(\a s + 1)}\Ga(\lbd) d\lbd \; =\; \frac{1}{\Gamma(-\a s)} \int_0^\infty \lbd^{-(\a s + 1)} (e^\lbd - \a\lbd^{\a -1}\Ea'(\lbd^\a)) d\lbd$$
but I could not carry the computations further with the sole help of Mittag-Leffler functions. Instead, one can rest upon the residue theorem in order to evaluate directly
$$\EE\lcr \Ua^s\rcr \; =\; \frac{-\sin \pi\a}{\pi}\int_0^\infty \frac{x^{s+1/\a}}{x^2 - 2x \cos\pi\a +1} dx.$$
Setting $\lbd = s +1/\a +1\in (0,2)$ and 
$$f(z)\; =\; \frac{(-z)^{\lbd -1}}{z^2 - 2z \cos\pi\a +1}$$
one sees from Formula VIII.(10.2.3) in \cite{D} that
$$\EE\lcr \Ua^s\rcr \; =\;-\lpa\frac{\sin \pi\a}{\sin\pi\lbd}\rpa\lpa {\rm Res}_{x_\a}f \, +\, {\rm Res}_{{\bar x}_\a}f\rpa$$
where $x_\a = e^{\ii(2-\a)\pi}$ and ${\bar x}_\a = e^{\ii(\a-2)\pi}$ are the two conjugate roots of $x^2 - 2x \cos\pi\a +1.$ After some simple computations, one gets
$$\EE\lcr \Ua^s\rcr \; =\;-\lpa\frac{\sin (\pi(\lbd -1)(\a-1))}{\sin\pi\lbd}\rpa\; =\; \frac{\sin ( \pi(\lbd -1)(\a-1))}{\sin\pi(\lbd-1)}\; =\; \frac{\sin (\pi (\a -1)(s+1/\a))}{\sin(\pi(s+1/\a))}$$
for every $s\in (-1-1/\a, 1-1/\a),$ which together with the previous computation for $\hta_1$ completes the proof (notice that $1- 1/\a < 1/\a$).
\fin

\begin{Rem} {\em The fractional moments of $\Ua$ and $\tau_1$ are finite over the same strip $(-1-1/\a, 1-1/\a),$ whereas those of $\hta_1$ are finite over the larger strip $(-\infty, 1/\a).$}
\end{Rem}

We will now derive series representations for the density $f^{}_{\tau_1}$ of $\tau_1$, inverting the fractional moments which were computed just before. Setting $M(z) = \EE[\tau_1^z]$ for 
every $z\in\CC$ such that ${\rm Re}(z) \in (-1-1/\a, 1-1/\a),$ one has 
$$M(\ii \lbd) \; =\; \EE[e^{\ii \lbd \sigma_1}], \quad\lbd\in\rl,$$
where $\sigma_1 = \log \tau_1$ is a real random variable with finite polynomial moments by Corollary 1, and which is absolutely continuous with continuous density $f^{}_{\sigma_1}(x) = e^x f^{}_{\tau_1}(e^x)$ over $\rl.$ The Fourier inversion formula for $f^{}_{\sigma_1}$ entails then
\begin{equation}
\label{Mellin}
f^{}_{\tau_1} (x)\; =\;\frac{1}{2\pi x} \int_{\rl} M(\ii t) x^{-\ii t} dt, \quad x\ge 0,
\end{equation}
and we will start from this formula to obtain our series representation. The one at zero is particularly simple:

\begin{Proni} For every $x \ge 0$ one has
$$f_{\tau_1} (x)\; =\; \sum_{n\ge 1}\frac{\a x^{1/\a +n-1}}{\Gamma(-\a n) \Gamma(1/\a +n)}\cdot$$

\end{Proni}

\noindent
{\em Proof}: Suppose first $x<1.$ We compute the integral in (\ref{Mellin}) with the help of the usual contour $\Gamma_R$ joining $-R$ to $R$ along the real axis and $R$ to $-R$ along a half-circle in the trigonometric orientation, centered at the origin. Because $x< 1,$ the integral along the half-circle is easily seen to converge to 0 when $R\to \infty,$ so that it remains to consider the singularities of $L(t) = M(\ii t) x^{-\ii t}$ inside the contour. By Corollary 1, the latter are located at $\ii(n+1/\a),\; n\ge 0.$ After standard computations, one finds that 
\begin{eqnarray*}
{\rm Res}_{\ii(n+1/\a)}L & = & -\ii\lpa \frac{\Gamma(2+\a n) \sin(n\pi \a) x^{n+1/\a}}{\pi\Gamma(1+1/\a +n)}\rpa\; =\; -\ii\lpa\frac{\a x^{n+1/\a}}{\Gamma(-\a n)\Gamma(1/\a +n)}\rpa
\end{eqnarray*}
for any $n\ge 0.$ By the residue theorem, this finishes the proof for $x \in (0, 1)$ (notice that  ${\rm Res}_{\ii /\a}L = 0$) and by analyticity, one can then extend the formula over the whole $\rl^+.$

\fin

The series representation at infinity is slightly more complicated than the one at zero, involving actually two series.
Up to some painless normalization and after some simplifications on the Gamma function, it had already been  computed by Peskir \cite{Pe} with an entirely different,  probabilistic argument involving a previous computation made for $f_{S_1}^{}$ in \cite{BDP}. 

 \begin{Prosan} For every $x > 0$ one has
 $$f_{\tau_1} (x)\; =\; \sum_{n\ge 1}\frac{-\a x^{1/\a -n-1}}{\Gamma(\a n) \Gamma(1/\a -n)}\; +\; \sum_{n\ge 1}\frac{x^{-n/\a-1}}{\Gamma(-n/\a) n!} \cdot$$
\end{Prosan}

\noindent
{\em Proof}: By analyticity, it suffices to consider the case $x >1.$ We compute the integral in (\ref{Mellin}) with a contour ${\bar \Gamma}_R$ which is the reflection of $\Gamma_R$ with respect to the real axis. Because $x> 1,$ the integral along the half-circle converges to 0 when $R\to \infty,$ so that one needs again to consider the singularities of $L$ inside ${\bar \Gamma}_R$, which are located at $\ii(1/\a- n)$ and $-\ii n/\a, \; n\ge 0.$ The term
$$\sum_{n\ge 1}\frac{-\a x^{1/\a -n-1}}{\Gamma(\a n) \Gamma(1/\a -n)}$$
follows then from the computations of Proposition 2 in replacing $n$ by $-n$ and taking into account the clockwise orientation of ${\bar \Gamma}_R.$ The last term follows from the simple computation
$${\rm Res}_{\ii(n+1/\a)}L \; = \; \frac{\ii x^{-n/\a}}{\Gamma(-n/\a) n!} \cdot$$

\fin

\begin{Rem}{ \em The two above representations are integrable term by term, which yields two series representations for the distribution function of $\tau_1:$ 
$$\pb[\tau_1\le x ]\; =\; \sum_{n\ge 1}\frac{\a x^{1/\a +n}}{\Gamma(-\a n) \Gamma(1+1/\a +n)}$$
and
$$\pb[\tau_1\ge x ]\; =\; \sum_{n\ge 1}\frac{\a x^{1/\a -n}}{\Gamma(\a n) \Gamma(1+1/\a -n)}\; -\; \sum_{n\ge 1}\frac{x^{-n/\a}}{\Gamma(1-n/\a) n!} \cdot$$
One can also differentiate term by term in order to get expansions for the successive derivatives of $f_{\tau_1}:$ for any $p\ge 1$ one has
$$f^{(p)}_{\tau_1} (x)\; =\; \sum_{n\ge 1}\frac{\a x^{1/\a +n-1-p}}{\Gamma(-\a n) \Gamma(1/\a +n-p)}$$
and
$$f^{(p)}_{\tau_1} (x)\; =\; \sum_{n\ge 1}\frac{-\a x^{1/\a -n-1-p}}{\Gamma(\a n) \Gamma(1/\a -n-p)}\; +\; \sum_{n\ge 1}\frac{x^{-n/\a-1-p}}{\Gamma(-n/\a -p) n!} \cdot$$
The first term of all these expansions give the behaviour of the function at zero respectively at infinity. For instance, one has
$$f_{\tau_1}(x)\; \sim\; \frac{\a x^{1/\a}}{\Gamma(-\a)\Gamma(1+1/\a)} \quad \mbox{and}\quad f'_{\tau_1}(x)\; \sim\; \frac{\a x^{1/\a-1}}{\Gamma(-\a)\Gamma(1/\a)}$$ 
as $x\to 0^+,$ whereas
$$f_{\tau_1}(x)\; \sim\; \frac{x^{1/\a-2}}{\Gamma(\a-1)\Gamma(1/\a)} \quad\mbox{and}\quad f'_{\tau_1}(x)\; \sim\; \frac{-\a x^{1/\a-3}}{\Gamma(\a)\Gamma(1/\a-2)}$$
as $x\to +\infty.$ Notice finally that since $\a \in (1,2),$ the function $f_{\tau_1}$ is ultimately completely monotone at infinity and $f'_{\tau_1}$ ultimately completely monotone at zero.}
\end{Rem}

\section{Complements on the first exit time and the running supremum}

In this section, we will derive similar series representations for the density $f^{}_{T_1}$ of the first exit time $T_1$. As mentioned before, this random variable is connected to the supremum process $\{S_t, \, t\ge 0\},$ so that we will also get series representations for $S_1.$ Our basic tool is an identity in law analogous to our main result, which was proved in our previous article \cite{Th1}, Theorem 3. This identity reads
\begin{equation}
\label{running}
T_1\elaw \Ta \times \hta_1
\end{equation}
where is an independent random variable with density given by
$$f_{\Ta}(t)\; =\;\frac{-(\sin \pi\a) (1+t^{1/\a})}{\pi\a(t^2 - 2t \cos \pi\a  +1)}$$
over $\rl^+.$ As before, we begin in evaluating the fractional moments of $T_1$:

\begin{ProV} For every $s\in (-1, 1-1/\a)$ one has
$$\EE\lcr T_1^s\rcr\; =\; \frac{\Gamma(1 + s) \sin (\pi/\a)}{\Gamma(1+s\a) \sin(\pi(s+1/\a))}\cdot$$
\end{ProV}
\noindent
{\em Proof.} For every $s\in (-1, 1-1/\a)$ we have
\begin{eqnarray*}
\EE\lcr \Ta^s\rcr & = & \frac{1}{\a}\lpa\EE\lcr \Ua^{s-1/\a}\rcr\; +\;\EE\lcr \Ua^s\rcr \rpa\\
& = & \frac{1}{\a}\lpa \frac{\sin(\pi s(\a -1))}{\sin \pi s} \; -\; \frac{\sin(\pi s(\a -1) -1/\a)}{\sin (\pi (s +1/\a))} \rpa\\
& = & \frac{\sin(\pi s \a)\sin(\pi/\a)}{\a\sin(\pi s) \sin (\pi (s +1/\a))}
\end{eqnarray*}
where the last equality follows from tedious trigonometric
transformations. From (\ref{running}), the previous computation made for $\EE
\lcr \hta_1^s\rcr$, and the complement formula for the Gamma function, we get
$$\EE\lcr T_1^s\rcr\; =\; \frac{\sin (\pi/\a)}{\sin(\pi(s+1/\a)}\lpa \frac{\Gamma(1 - s\a) \sin (\pi s \a))}{\a \Gamma(1- s) \sin(\pi s)}\rpa\; =\; \frac{\Gamma(1 + s) \sin (\pi/\a)}{\Gamma(1+s\a) \sin(\pi(s+1/\a))}$$
as desired.
\fin

\begin{Rem}{\em Similarly as above, the fractional moments of $\Ta$ and $\tau_1$ are finite over the same strip $(-1, 1-1/\a),$ whereas those of $\hta_1$ are finite over the larger strip $(-\infty, 1/\a).$}
\end{Rem}

We notice in passing that this computation allows to show that the law of the independent quotient $T_1/\hT_1$ is ${\rm Pareto} (1-1/\alpha),$ a fact which had been proved in \cite{Do1} in the general context of stable processes in duality:

\begin{Cogo} The density of the independent quotient $T_1/\hT_1$ is given by 
$$\frac{\sin(\pi/\alpha)}{\pi t^{1-1/\alpha}(1+t)}$$
over $\rl^+.$
\end{Cogo}

\noindent
{\em Proof}: By Mellin inversion, it is enough to identify the fractional moments. Proposition 4 and the classical computation for the moments of $\hT_1$ entail
$$\EE\lcr \frac{T_1^s}{\hT_1^s}\rcr\; =\; \frac{\sin (\pi/\a)}{\sin(\pi(s+1/\a))}\; =\; \int_0^\infty \frac{\sin(\pi/\alpha) t^s}{\pi t^{1-1/\alpha}(1+t)} dt$$
for every $s\in (-1/\a, 1- 1/\a),$ where the second equality follows from Formula VIII.(10.2.4) in \cite{D}.

\fin

We now come back to series representations for $f_{T_1}^{}$, which are obtained from Proposition 4 and Mellin inversion similarly as in Propositions 2 and 3. Since we choose exactly the same contours, we will leave all the details to the reader. Notice that contrary to $\tau_1,$ the two-series representation for $T_1$ is the one at zero and that the first series therein defines an entire function.

\begin{Corroku} One has both representations
$$f^{}_{T_1}(x)\; =\; \sum_{n\ge 1}\frac{x^{1/\a - n -1}}{\a\Gamma(\a n -1)\Gamma(1+1/\a -
  n)}, \quad x > 0$$
and
$$f^{}_{T_1}(x)\; =\; \sum_{n\ge 1}\frac{x^{n -1}}{\a\Gamma(-\a n)}\; +\; \sum_{n\ge 1}\frac{x^{1/\a + n -1}}{\Gamma(-\a n)\Gamma(1/\a + n)}, \quad x \ge 0.$$
\end{Corroku}

Thanks to the aforementioned identity $T_1 \elaw S_1^{-\a},$ changing the variable
gives corresponding series representations for the density $f^{}_{S_1} = \a x^{-(\a +1)}f^{}_{T_1}(x^{-\a}).$
The first one had been obtained originally in \cite{BDP} after solving some
fractional integral equation of the Abel type. See also \cite{P} for a new
proof involving Wright's hypergeometric function. We stress that up to the
Wiener-Hopf factorisation which is a fundamental tool for this type of
questions, overall our argument to get this series only makes use of undergraduate mathematics. Besides, this simple method allows to get the series representation at infinity.

\begin{Coshichi} One has both representations
$$f^{}_{S_1}(x)\; =\; \sum_{n\ge 1}\frac{x^{\a n -2}}{\Gamma(\a n -1)\Gamma(1+1/\a -
  n)}, \quad x\ge 0$$
and
$$f^{}_{S_1}(x)\; =\; \sum_{n\ge 1}\frac{x^{-(\a n +1)}}{\Gamma(-\a n )}\; +\; \sum_{n\ge 1}\frac{\a x^{-(\a n +2)}}{\Gamma(-\a n)\Gamma(1/\a + n)}, \quad x > 0.$$
\end{Coshichi}

\begin{Rem} {\em Notice that one really needs the two series to get the behaviour of $f^{}_{T_1}$ and its derivatives at zero. For instance one has
$$f^{}_{T_1}(x)\; \sim\; \frac{1}{\a\Gamma(-\a)} \quad \mbox{and}\quad f'_{T_1}(x)\; \sim\; \frac{x^{1/\a-1}}{\Gamma(-\a)\Gamma(1/\a)}\;\sim\; \frac{1}{\a} f'_{\tau_1}(x)$$ 
as $x\to 0^+.$ On the other hand, only one series is important to get the behaviour of $f^{}_{S_1}$ and its derivatives at infinity: one has
$$f^{}_{S_1}(x)\; \sim\; \frac{x^{-(\a +1)}}{\Gamma(-\a)} \quad \mbox{and}\quad f^{(p)}_{S_1}(x)\; \sim\; \frac{x^{-(\a+2+p)}}{\Gamma(-(\a +1+p))}$$ 
as $x\to +\infty$ for every $p\ge 1$ (this is in accordance with the general asymptotics derived in \cite{DS}, see the final remark therein).}
\end{Rem}

\section{Some remarks on unimodality}

Recall that a real random variable $X$ is said to be unimodal if there exists
$a\in\rl$ such that the functions $\pb[X\le x]$ and $\pb[X>x]$ are convex
respectively in $(-\infty,a)$ and $(a,+\infty).$ If $X$ has a density $f_X^{}$, this means that $f^{}_X$ increases on $(-\infty,a]$ and decreases on $[a,+\infty).$ We refer e.g. to Section 52 in \cite{S} for more on this topic. Differentiating the density of $\Ua:$ 
$$f_{\Ua}'(t)\; =\; \frac{t^{1/\a-1}}{\a(t^2 -2t\cos\pi\a +1)^2}((1-2\a)t^2 + 2t(\a -1) \cos\pi\a +1),$$
we find that there is a unique positive root, so that the random variable $\Ua$ is unimodal. Since the same property is known to hold for $\hta_1$ - see e.g. Theorem 53.1 in \cite{S} for a much stronger result, in view of our main result it is natural to ask whether the variable $\tau_1$ is unimodal as well. This seems also plausible from the behaviour of $f_{\tau_1}$ and $f_{\tau_1}'$ at zero. In general, the product of two positive independent unimodal random variables is not necessarily unimodal. However the notion of multiplicative strong unimodality which had been introduced in \cite{CT}, will allow us to give a quick positive answer in half of the cases:

\begin{Proachi} For every $\a \le 3/2,$ the random variable $\tau_1$ is unimodal.
\end{Proachi}
\noindent
{\em Proof}: Recalling that $\hta_1$ is unimodal for every $\a\in (1, 2),$
from (\ref{Main}) and Theorem 3.6 in \cite{CT} is is enough to show that the function $t\mapsto f_{\Ua}(e^t)$ is log-concave over $\rl,$ in other words that the function $g_\a (t) = \log (e^t - 2\cos\pi\a + e^{-t})$ is convex over $\rl.$ Differentiating twice yields 
$$g_\a'' (t)\; =\;\frac{4(1 -\cos \pi\a \cosh t)}{(e^t - 2 \cos \pi\a + e^{-t})^2}$$
and we see that $g_\a$ is convex over $\rl$ as soon as $\a\le 3/2.$

\fin

When $\a \in (3/2, 2)$ the above proof shows that $\Ua$ is not multiplicative
strongly unimodal (MSU) in the terminology of \cite{CT}. In order to get the
unimodality of $\tau_1$ when $\a \in (3/2, 2),$ one could of course be tempted
to obtain the MSU property for positive stable laws with index $\beta > 2/3.$
Unfortunately, we showed in \cite{Th2} that this property does not hold true
as soon as $\beta > 1/2.$ Nevertheless, in view of the unimodality of $\hta_1$
and that of $\tau_1$ for $\a = 2$ - see also R\"osler's result for general
diffusions on the line \cite{R}, we conjecture that $\tau_1$ is unimodal for every value of $\a.$ At first sight the problem does not seem easy because of the unstability of the unimodality property under multiplicative convolution. See however our final remark. Differentiating now the density of $\Ta:$ 
$$f_{\Ta}'(t)\; =\; \frac{t^{1/\a-1}}{\a(t^2 -2t\cos\pi\a +1)^2}((1-2\a)t^2 - 2\a t^{2-1/\a} + 2t(\a -1) \cos\pi\a + 2\a t^{1-1\a}\cos \pi\a  +1),$$
we find after some simple analysis that there is also a unique positive root, so that the random variable $\Ta$ is unimodal. But contrary to $\Ua,$ one can check that it is never MSU, so that the above simple argument cannot be applied. Nevertheless we have the

\begin{Prokyu} For every $\a \le 3/2,$ the random variable $T_1$ is unimodal.
\end{Prokyu}
\noindent
{\em Proof}: From (\ref{running}) and the so-called Kanter representation for positive stable laws - see Corollary 4.1 in \cite{K}, we have the identity
$$T_1\;\elaw\; \Ta\times L^{1-\a} \times b_\a(U)$$
where $L$ is a standard exponential variable, $U$ an independent uniform
variable over $[0,\pi],$ and $b_\a(u) = (\sin((1-1/\a)u)/\sin
(u))^{\a-1}\sin(u/\a)/\sin(u)$ is a strictly increasing function from
$(0,\pi)$ to $\rl^+$. Besides, from the beginning of the proof of Theorem 4.1
in \cite{K}, we know that the logarithmic derivative of $b_\a$ is positive and
strictly increasing, so that the same holds for $b_\a'$ itself, because $b_\a$
is also positive and strictly increasing. By Lemma 4.2 in \cite{K} we conclude
that $b_\a(U)$ is unimodal: again, from (\ref{running}) and Theorem 3.6 in \cite{CT}, it suffices to show that the density $t\mapsto f_{\Ta\times L^{1-\a}}^{}(e^t)$ is log-concave over $\rl.$ We compute 
$$f_{\Ta\times L^{1-\a}}^{} (x)\; =\; \frac{-(\sin \pi\a)x^{\frac{\a}{1-\a}}}{\pi\a(\a - 1)} \int_0^\infty \exp-\lpa x/u\rpa^{\frac{1}{1-\a}}\frac{u^{\frac{1}{\a -1}}(1 + u^{1/\a})}{u^2 -2u\cos \pi\a +1} du, \quad x \ge 0.$$
It is easily seen that $u\mapsto 1+ u^{1/\a}$ is log-concave over $\rl^+$ and $(t,u)\mapsto e^t/u$ log-convex, hence convex, over $\rl\times\rl^+$. In particular, the function 
$$(t,u)\;\mapsto\;\exp-\lpa e^t/u\rpa^{\frac{1}{1-\a}}(1 + u^{1/\a})$$ 
is log-concave over $\rl\times\rl^+$ and by Pr\'ekopa's theorem, it is enough to show that the function
$u\mapsto u^{\frac{1}{\a -1}}/(u^2 -2u\cos \pi\a +1)$
is log-concave over $\rl^+.$ We compute its logarithmic derivative, which equals
$$\frac{(2\a -3)u^4 +4u^3(2-\a)\cos\pi \a - (4(2-\a)\cos^2\!\pi \a + 2\a) u^2 + 4u\cos\pi \a - 1}{(\a -1)u^2(u^2 -2u\cos \pi\a +1)^2}$$ 
and we see that it is negative over $\rl^+$ as soon as $\a \le 3/2.$

\fin

\begin{Frem}{\em The above proof does not work to obtain the unimodality of
    $T_1$ when $\a > 3/2,$ and we do not know as yet how to tackle this
    situation. However, we can use the same argument to obtain the unimodality
    of $\tau_1$ when $3/2 <\a \le 1 + 1/\sqrt{2} <2$. Details, which are
    technical, can be obtained upon request. All these questions will be the
    matter of future research.}
\end{Frem}

\bigskip

\noindent
{\bf Acknowledgements.} Part of this work was done during a sunny stay at the
University of Tokyo and I am very grateful to Nakahiro Yoshida for his hospitality. Ce travail a
aussi b\'en\'efici\'e d'une aide de l'Agence Nationale de la Recherche
portant la r\'ef\'erence ANR-09-BLAN-0084-01 (projet {\em Autosimilarit\'e}).

\end{document}